\documentclass{amsart}

\usepackage{epsfig}
\usepackage{graphics}
\usepackage{amsmath,amsfonts}

\newtheorem{theorem}{Theorem}[section]
\newtheorem{lemma}[theorem]{Lemma}
\newtheorem{prop}[theorem]{Proposition}

\theoremstyle{definition}
\newtheorem{definition}[theorem]{Definition}

\theoremstyle{remark}
\newtheorem{remark}[theorem]{Remark}

\numberwithin{equation}{section}

\newcommand{\ba}{${\scriptstyle /}$}

\newcommand{\Cs}{C \kern-2.2mm \raise.4ex\hbox {\ba} \kern.2mm}
\newcommand{\Ns}{N \kern-3.4mm \raise.35ex\hbox {\ba} \kern.4mm\hskip 1.3mm}
\newcommand{\bas}{${\scriptscriptstyle /}$}
\newcommand\Css{C \kern-1.9mm \raise.3ex\hbox {\bas} \kern.8mm}

\begin{document}

\title{Multi-point Taylor Expansions of Analytic Functions}

\author{Jos\'e L. L\'opez}
\address{Departamento de Mat\'ematica e Inform\'atica,
Universidad P\'ublica de Navarra, 31006-Pamplona,
Spain}
\email{jl.lopez@unavarra.es}
\thanks{The first author wants to thanks the saving bank {\it Caja Rural
de Navarra} for its financial support. He also acknowledges
the scientific and financial support of CWI in Amsterdam.}

\author{Nico M. Temme}
\address{CWI, P.O. Box 94079, 1090 GB Amsterdam, The Netherlands}
\email{nicot@cwi.nl}
\thanks{The authors thank the referee for the comments on the
first version of the paper.}


\date{October 1, 2002 and, in revised form, January 20, 2004.}



\begin{abstract}
Taylor expansions of analytic functions are considered with respect to several
points, allowing confluence of any of them.
Cauchy-type formulas are given for coefficients and remainders in
the expansions, and the regions of convergence are indicated. It is
explained how these expansions can be used in deriving  uniform asymptotic
expansions of integrals. The method is also used for obtaining  Laurent
expansions in several points as well as Taylor-Laurent expansions.
\end{abstract}

\maketitle



\section{Introduction}
In deriving uniform asymptotic expansions of a certain class of integrals
one encounters the problem of expanding a function, that is analytic in some
domain $\Omega$ of the complex plane, in several points. The first mention
of the use of such expansions in asymptotics is given in \cite{Chester},
where Airy-type expansions are derived for integrals having two
nearby (or coalescing) saddle points. This reference does not give further
details about two-point Taylor expansions, because the coefficients in the
Airy-type asymptotic expansion are derived in a different way. Other mentions
of the use of such expansions in asymptotics is given in \cite{Raimundas}
and \cite{nicoi}.
In \cite{Raimundas}, two-point Taylor expansions are used
with applications to Airy-type expansions of parabolic cylinder functions.
In \cite{nicoi} we used two-point Taylor expansions to derive
convergent expansions of Charlier, Laguerre and Jacobi polynomials in terms
of Gamma, Hermite and Chebyshev polynomials respectively.

To demonstrate an application in asymptotics of multi-Taylor expansions
we may consider contour integrals of the form
\begin{equation}
\label{intro}
I(\lambda;{\bf\alpha})=\int_{C} g(z)e^{-\lambda f(z,{\bf\alpha})}\,dz,
\end{equation}
where ${\bf\alpha}$ is a vector of parameters,
${\bf\alpha}=(\alpha_1,\ldots,\alpha_2)$ and the phase function
$f(z,{\bf\alpha})$
has $m$ saddle points $z_1$, $z_2$,...,$z_m$. The asymptotic behaviour of these
integrals for large values of $\lambda$ is determined by
the saddle-point structure of the phase function
[\cite{wong}, chap. 7, sec. 6]. One method for obtaining an
asymptotic expansion of this integral for large values of $\lambda$
is based on expanding $g(z)$ at the saddle points of the phase function,
\begin{equation}
g(z)=\sum_{n=0}^\infty[a_0+a_1z+\ldots+a_{m-1}z^{m-1}]
(z-z_1)^n(z-z_2)^n\cdots(z-z_m)^n
\end{equation}
and substitute this expansion into (\ref{intro}). When interchanging summation
and integration, the result is a formal expansion in $m$ series in
terms of functions related with the functions
\begin{equation}
F_{n,k}(\lambda;{\bf\alpha})\equiv\int_{C}
z^k(z-z_1)^n(z-z_2)^n\cdots(z-z_m)^ne^{-\lambda f(z,{\bf\alpha})}\,dz,
\hskip 5mm k=1,2,\ldots,m-1.
\end{equation}
In \cite{Raimundas}, these functions $F_{n,k}(\lambda;{\bf\alpha})$
are the Airy functions,
whereas in \cite{nicoi} these functions are the Gamma function,
or the Hermite or Chebyshev polynomials.

In a future paper we will use multi-point Taylor expansions in
the asymptotic analysis of integrals arising in diffraction theory, such as
the Bessel function integral (see \cite{Janssen} and \cite{Kirk})
\begin{equation}
\label{besel}
J(x,y)=\int_0^\infty t J_0(yt) e^{i(t^4+xt^2)}\,dt,
\end{equation}
which is related to the Pearcey function
\begin{equation}
\label{pearcy}
\int_{-\infty}^{\infty} e^{i({1\over 4}t^4+{1\over 2}xt^2+yt)}dt.
\end{equation}

The Taylor-Laurent expansions
will be used to study integrals with two saddle points and a pole of
the integrand. Other applications in asymptotics include the
study of Hermite-Pad\'e
approximations to the exponential function; in \cite{Driver} integrals are
considered
with three saddle points.


In a recent paper \cite{nicoii} we have introduced the theory of
two-point Taylor expansions, two-point Laurent expansions and
two-point Taylor-Laurent expansions. The purpose of the
present paper is to generalize that theory from 2 to $m$ points,
$m\ge 2$.
We give details on the region of convergence and on
representations in terms of Cauchy-type integrals of the coefficients and
the remainders of the expansions.
Earlier information on this type of expansions is given
in \cite{walsh}, Chapters 3 and 8. The theory of
several-point Taylor expansions was already formulated in Chapter 3 of
Walsh's book, although in a different setting.
Chapter 8  of \cite{walsh} presents also a theory of rational approximation of
analytic functions, but is different from the theory of multi-point Laurent
and Taylor-Laurent expansions presented here. Whereas the
multi-point polynomial approximation of Chapter 3 may be reformulated as a
multi-point Taylor approximation, the rational approximation of Chapter 8
can not be written as a multi-point Laurent or Taylor-Laurent approximation.
For more details, see Section 5.

\section{Multy-point Taylor expansions}
\noindent
We consider the Taylor expansion of an analytic function $f(z)$ in several
points and give information on
the coefficients and the remainder in the expansion. In what follows
empty sums and derivatives of negative order must be understood as
zero and empty products as one. We
will deal with the following set of points:

\begin{definition}
\label{def1}
     We define the set
\begin{equation}
S\equiv\lbrace z_1,z_1,\ldots,z_1;
z_2,z_2,\ldots,z_2,\ldots;z_p,z_p,\ldots,z_p\rbrace
\end{equation}
of $m$ points consisting on $p$ different points $z_1$, $z_2$,..., $z_p$
($z_i\ne z_j$ if $i\ne j$),
each $z_j$ repeated $m_j$ times: $m_1+m_2+\ldots+m_p=m$.
\end{definition}

For clearity in the exposition, we first introduce the multi-point Taylor
expansion for $m$ different points $z_1$, $z_2$, $\ldots z_m$ ($m=p$, $m_j=1$)
in Theorem \ref{theo1}. In Theorem \ref{theo2} we assume that the points
$z_1$, $z_2$, $\ldots z_m$ may coalesce.
We will need the following elementary lemma.

\begin{lemma}
\label{lema1}
     Given $z$, $w$, $\in\Cs$, take $m$
different points $z_1$, $z_2$,$\ldots$, $z_m$ in $\Cs$ and define
\begin{equation}
     \label{defache}
H_m(w,z;z_1,\ldots,z_m)\equiv{\prod_{k=1}^m(w-z_k)
-\prod_{k=1}^m(z-z_k)\over w-z}.
\end{equation}
Then
\begin{equation}
     \label{ache}
H_m(w,z;z_1,\ldots,z_m)=\sum_{j=1}^m{\prod_{k=1,k\ne j}^m(w-z_k)
\prod_{k=1,k\ne j}^m(z-z_k)\over \prod_{k=1,k\ne j}^m(z_j-z_k)}.
\end{equation}
\end{lemma}

\begin{proof}
     The numerator of $H_m(w,z;z_1,\ldots,z_m)$ is a polynomial of
degree $m$ in the variable $w$ that vanishes at $w=z$. Therefore,
$H_m(w,z;z_1,\ldots,z_m)$ is a polynomial of degree $m-1$ in the variable $w$.
Let $P_m(w,z;z_1,\ldots,z_m)$ denote the function at the right-hand side of
(\ref{ache}), which is also
a polynomial of degree $m-1$ in the variable $w$. Moreover,
\begin{equation}
H_m(z_s,z;z_1,\ldots,z_m)=\prod_{k=1,k\ne
s}^m(z-z_k)=P_m(z_s,z;z_1,\ldots,z_m)
\end{equation}
for $s=1$, $2$,$\ldots m$. Hence,
\begin{equation}
H_m(w,z;z_1,\ldots,z_m)=P_m(w,z;z_1,\ldots,z_m).
\end{equation}
\end{proof}

\begin{theorem}
\label{theo1}
Let $f(z)$ be an analytic function on an open set
$\Omega\subset C$ and $S\subset\Omega$ with $S$ consisting
on $m$ different points ($m=p$).
Then, $f(z)$ has the multi-point Taylor expansion
\begin{equation}
\label{expani}
f(z)=\sum_{n=0}^{N-1}q_{n,m}(z)\prod_{k=1}^m(z-z_k)^{n}+r_N(z),
\end{equation}
where $q_{n,m}(z)$ is the polynomial of degree $m-1$
\begin{equation}
\label{poly}
q_{n,m}(z)\equiv\sum_{j=1}^m a_{n,j}
{\prod_{k=1,k\ne j}^m(z-z_k)\over \prod_{k=1,k\ne j}^m(z_j-z_k)}
\end{equation}
and the coefficients $a_{n,j}$ of this polynomial are given by the
Cauchy integral
\begin{equation}
\label{coefi}
a_{n,j}\equiv{1\over 2\pi i}
\int_C{f(w)\,dw\over (w-z_j)\prod_{k=1}^m(w-z_k)^{n}}.
\end{equation}
The remainder term $r_N(z)$ is given by the Cauchy integral
\begin{equation}
     \label{remi}
r_N(z)\equiv{1\over 2\pi i}\int_C{f(w)\,dw\over (w-z)
\left[\prod_{k=1}^m(w-z_k)\right]^N}
\left[\prod_{k=1}^m(z-z_k)\right]^N.
\end{equation}
The contour of integration $C$ is a simple closed loop which
encircles the points $z_1$, $z_2$,...,$z_m$ (for $a_{n,j}$)
and $z$, $z_1$, $z_2$,...,$z_m$ (for $r_N(z)$) in
the counterclockwise direction and is contained in $\Omega$ (see
Figure 1 (a)).

The expansion (\ref{expani}) is convergent for $z\in O_m$, where:
\begin{equation}
\label{ere}
O_m\equiv\lbrace z\in\Omega, \hskip 2mm
\prod_{k=1}^m\vert z-z_k\vert<r\rbrace, \hskip 2cm
r\equiv {\rm Inf}_{w\in\Css\setminus\Omega}\left\lbrace
\prod_{k=1}^m\vert w-z_k\vert\right\rbrace.
\end{equation}
That is, (\ref{expani}) is convergent for $z$ inside the lemniscate
$\prod_{k=1}^m\vert z-z_k\vert=r$ (see
Figure 2; if $m=1$ this domain is a disk;
if $m=2$ this domain is bounded by a Cassini oval).
In particular, if $f(z)$ is an entire function $(\Omega=\Cs$), then the
expansion (\ref{expani}) converges $\forall$ $z\in\Cs$.
\end{theorem}

\begin{proof}
By Cauchy's theorem,
\begin{equation}
\label{cauchy}
f(z)={1\over 2\pi i}\int_C{f(w)\,dw\over w-z},
\end{equation}
where $C$ is the contour defined above (Figure 1 (a)). We write
\begin{equation}
\label{lequal}
{1\over w-z}={H_m(w,z;z_1,\ldots,z_m)\over \prod_{k=1}^m(w-z_k)}
{1\over 1-u},
\end{equation}
where $H_m(w,z;z_1,\ldots,z_m)$ is given in (\ref{defache}) and
\begin{equation}
     \label{ui}
u\equiv {\prod_{k=1}^m(z-z_k)\over \prod_{k=1}^m(w-z_k)}.
\end{equation}
Now we use Lemma \ref{lema1} and introduce the right hand side of
(\ref{ache}) and the expansion
\begin{equation}
     \label{expanu}
{1\over 1-u}=\sum_{n=0}^{N-1}u^n+{u^N\over 1-u}
\end{equation}
in (\ref{lequal}) and this in (\ref{cauchy}). After straightforward 
calculations we
obtain formulas (\ref{expani})-(\ref{remi}).

For any $z\in O_{m}$, we can take a contour $C$ in $\Omega$ such that
\begin{equation}
\prod_{k=1}^m\vert z-z_k\vert<\prod_{k=1}^m\vert w-z_k\vert,
\quad \forall w\in\Cs
\end{equation}
(see Figure 1 (b)). On this contour $\vert f(w)\vert$ is bounded by
some constant $C$: $\vert
f(w)\vert\le C$. Introducing these two bounds in (\ref{remi}) we see that
$\lim_{N\to\infty}r_N(z)=0$ and the proof follows.
\end{proof}

\begin{figure}[tb]
     \begin{center}\includegraphics[width=8cm]{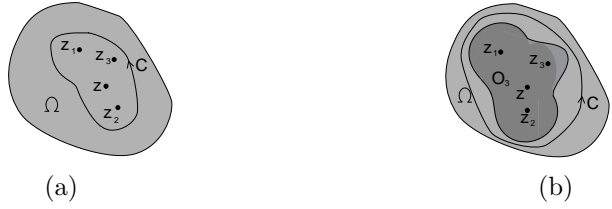}\end{center}
\centerline{(a) \hskip 6cm (b)}
\caption{The case $m=3$. (a) Contour $C$ in the integrals
(\ref{coefi}) and (\ref{remi}). (b) For
$z\in O_{m}$, we can take a contour $C$ in $\Omega$ which contains
$O_{m}$ inside and therefore,
$\prod_{k=1}^m\vert z-z_k\vert<\prod_{k=1}^m\vert w-z_k\vert$
$\forall$ $w\in\Cs$.}
\label{Figure 1}
\end{figure}

\begin{figure}[tb]
     \begin{center}\includegraphics[width=12cm]{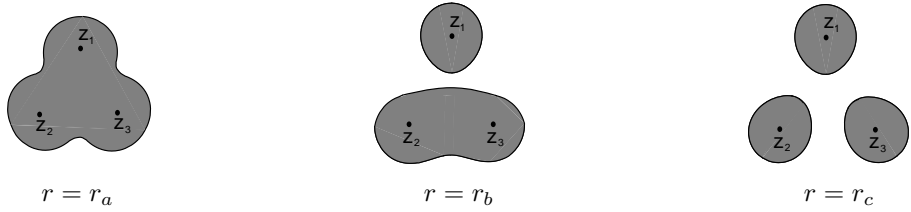}\end{center}
\centerline{$r=r_a$ \hskip 4cm $r=r_b$ \hskip 4cm $r=r_c$}
\caption{Shape of the "lemniscate domain" $O_m$ for $m=3$. It
depends on the size of the parameter $r$ defined in (\ref{ere}). In 
these pictures
$\vert z_2-z_3\vert<\vert z_1-z_3\vert$, $\vert z_1-z_2\vert$
and $r_a>r_b>r_c$.}
\label{Figure 2}
\end{figure}

We need the following lemma to consider the case of coalescing points in the
set~$S$.

\begin{lemma}
\label{lema2}
     Given $z,w\in\Cs$, take $m$ different points $z_1$, $z_2$,
..., $z_m$ in $C$, all different from $w$ too. Let those $m$ points
to coalesce at $z_m$, say. Then
\begin{equation}
\lim_{z_1,z_2,..,z_{m-1}\to z_m}\sum_{j=1}^m{\prod_{k=1,k\ne j}^m(z-z_k)
\over (w-z_j)\prod_{k=1,k\ne j}^m(z_j-z_k)}=
\sum_{j=0}^{m-1}{(z-z_m)^j\over (w-z_m)^{j+1}}.
\end{equation}
\end{lemma}

\begin{proof}
     We first note that the identity
\begin{equation}
\sum_{j=1}^n\prod_{l=1}^{j-1}(z_n-z_l)\prod_{l=j+1}^n(z_1-z_l)=0
\end{equation}
holds for any set of points $z_1$, $z_2$,..., $z_n$, $n>1$. It may be
cheked in the following way: we take the two first terms of the sum,
which gives
\begin{equation}
(z_n-z_2)(z_1-z_3)(z_1-z_4)\cdot\cdot\cdot(z_1-z_{n-1})(z_1-z_n).
\end{equation}
Next we add
to this the third term of the sum, which gives
\begin{equation}
(z_n-z_2)(z_n-z_3)(z_1-z_4)\cdot\cdot\cdot(z_1-z_{n-1})(z_1-z_n).
\end{equation}
We continue
this process untill we add the $n-1$-th term of the sum, obtaining
\begin{equation}
(z_n-z_2)(z_n-z_3)(z_n-z_4)\cdot\cdot\cdot(z_n-z_{n-1})(z_1-z_n).
\end{equation}
But this is just the last term of the sum with opposite sign.

Using the above identity we have
\begin{equation}
\prod_{l=1}^{k-1}(z_s-z_l)\left\lbrack
\sum_{j=k}^s\prod_{l=k}^{j-1}(z_s-z_l)\prod_{l=j+1}^s(z_k-z_l)
\right\rbrack \prod_{l=s+1}^{m}(z_k-z_l)=0
\end{equation}
for any $s=1,2,3,...,m$ and any $k=1,2,3,...,m$ with $k\ne s$. Then
\begin{equation}
     \label{cero}
\sum_{j=k}^s\prod_{l=1}^{j-1}(z_s-z_l)\prod_{l=j+1}^m(z_k-z_l)=
\sum_{j=1}^m\prod_{l=1}^{j-1}(z_s-z_l)\prod_{l=j+1}^m(z_k-z_l)=0
\end{equation}
for any $s,k=1,2,3,...,m$ with $k\ne s$.

Now, for every $s=1,2,3,...,m$, we define the
following polynomials of degree $m-1$ in the variable $z$:
\begin{equation}
R_s(z)\equiv\prod_{l=1,l\ne s}^m(z-z_l),
\hskip 2cm
S_s(z)\equiv\sum_{j=1}^m\prod_{l=1}^{j-1}(z-z_l)\prod_{l=j+1}^m(z_s-z_l).
\end{equation}
The zeros of $R_s(z)$ are $z_k$ for $k=1,2,3,...,m$, $k\ne s$ and from
(\ref{cero}),
$S_s(z_k)=0$ for $k=1,2,3,...,m$, $k\ne s$. Moreover, the leading
coeficient of $R_s(z)$
and of $S_s(z)$ coincide. Therefore, $R_s(z)=S_s(z)$ for $s=1,2,3,...,m$.

Finally, define the following polynomials of degree $m-1$ in the variable $w$:
\begin{equation}
P(w,z)\equiv\sum_{j=1}^m{\prod_{k=1,k\ne j}^m(z-z_k)\prod_{k=1,k\ne
j}^m(w-z_k)\over
\prod_{k=1,k\ne j}^m(z_j-z_k)},
\end{equation}
\begin{equation}
Q(w,z)\equiv\sum_{j=1}^m\prod_{k=1}^{j-1}(z-z_k)\prod_{k=j+1}^m(w-z_k).
\end{equation}
For every $s=1,2,3,...,m$ we have $P(z_s,z)=R_s(z)$ and $Q(z_s,z)=S_s(z)$. But
$R_s(z)=S_s(z)$ and therefore $P(w,z)=Q(w,z)$. Then,
\begin{equation}
\begin{split}
\sum_{j=1}^m{\prod_{k=1,k\ne j}^m(z-z_k)\over(w-z_j)\prod_{k=1,k\ne
j}^m(z_j-z_k)}= &
{P(w,z)\over\prod_{k=1}^m(w-z_k)}= \\ & {Q(w,z)\over\prod_{k=1}^m(w-z_k)}=
\sum_{j=1}^m{\prod_{k=1}^{j-1}(z-z_k)\over\prod_{k=1}^j(w-z_k)}.
\end{split}
\end{equation}
Taking the limits $z_1$, $z_2$,.., $z_{m-1}$ $\to$ $z_m$ in the left- and
right-hand sides of these equalities, we obtain the desired result.
\end{proof}

\begin{theorem}
\label{theo2}
     Let $f(z)$ be an analytic function on an open set
$\Omega\subset\Cs$ and $S\subset\Omega$.
Then, $f(z)$ has the multi-point Taylor expansion
\begin{equation}
     \label{expanib}
f(z)=\sum_{n=0}^{N-1}q_{n,m}(z)\prod_{k=1}^p(z-z_k)^{nm_k}+r_N(z),
\end{equation}
where $q_{n,m}(z)$ is the polynomial of degree $m-1$
\begin{equation}
     \label{polyb}
q_{n,m}(z)\equiv\sum_{j=1}^p
{\prod_{k=1,k\ne j}^p(z-z_k)^{m_k}\over \prod_{k=1,k\ne j}^p(z_j-z_k)^{m_k}}
\sum_{l=0}^{m_j-1}a_{n,j,l}(z-z_j)^l
\end{equation}
and the coefficients $a_{n,j,l}$ of this polynomial are given by the
Cauchy integral
\begin{equation}
     \label{coefib}
a_{n,j,l}\equiv{1\over 2\pi i}
\int_C{f(w)\,dw\over (w-z_j)^{l+1}\prod_{k=1}^p(w-z_k)^{nm_k}}.
\end{equation}
The remainder term $r_N(z)$ is given by the Cauchy integral
\begin{equation}
     \label{remib}
r_N(z)\equiv{1\over 2\pi i}\int_C{f(w)\,dw\over (w-z)
\left[\prod_{k=1}^p(w-z_k)^{m_k}\right]^N}
\left[\prod_{k=1}^p(z-z_k)^{m_k}\right]^N.
\end{equation}
The contour of integration $C$ is a simple closed loop which
encircles the points $z_1$, $z_2$,...,$z_p$ (for $a_{n,j,l}$)
and $z$, $z_1$, $z_2$,...,$z_p$ (for $r_N(z)$) in
the counterclockwise direction and is contained in $\Omega$ (see
Figure 1 (a)).

The expansion (\ref{expanib}) is convergent for $z\in O_p$:
\begin{equation}
O_p\equiv\lbrace z\in\Omega, \hskip 2mm
\prod_{k=1}^p\vert z-z_k\vert^{m_k}<r\rbrace, \hskip 2cm
r\equiv {\rm Inf}_{w\in\Css\setminus\Omega}\left\lbrace
\prod_{k=1}^p\vert w-z_k\vert^{m_k}\right\rbrace,
\end{equation}
that is, inside the lemniscate $\prod_{k=1}^p\vert z-z_k\vert^{m_k}=r$.
In particular, if $f(z)$ is an entire function $(\Omega=\Cs$), then the
expansion (\ref{expanib}) converges $\forall$ $z\in\Cs$.
\end{theorem}

\begin{proof}
     If all the points in $S$ are different,
we have from (\ref{poly}) and (\ref{coefi})
\begin{equation}
     \label{cuco}
q_{n,m}(z)={1\over 2\pi i}
\int_C{f(w)\,dw\over \prod_{k=1}^m(w-z_k)^n}
\sum_{j=1}^m{\prod_{k=1, k\ne j}^m(z-z_k)
\over (w-z_j)\prod_{k=1, k\ne j}^m(z_j-z_k)}.
\end{equation}
This last sum may be also decomposed in the form
\begin{equation}
\begin{split} &
\sum_{j=1}^{m_1}{\prod_{k=1, k\ne j}^m(z-z_k)
\over (w-z_j)\prod_{k=1, k\ne j}^m(z_j-z_k)}+
\sum_{j=m_1+1}^{m_2}{\prod_{k=1, k\ne j}^m(z-z_k)
\over (w-z_j)\prod_{k=1, k\ne j}^m(z_j-z_k)}+\ldots \\ &
\sum_{j=m_{p-1}+1}^m{\prod_{k=1, k\ne j}^m(z-z_k)
\over (w-z_j)\prod_{k=1, k\ne j}^m(z_j-z_k)}.
\end{split}
\end{equation}
Now let the first $m_1$ points to coalesce to $z_1$, the second $m_2$
points to coalesce to $z_2$, and so on, and apply Lemma \ref{lema2} to every
one of the $p$
sums above to obtain (\ref{expanib}), (\ref{polyb}) and 
(\ref{coefib}). Equation
(\ref{remib}) follows from (\ref{remi}). The proof of the convergence
of (\ref{expanib}) in the region $O_{p}$ is a straightforward generalization
of the correponding proof in Theorem \ref{theo1}.
\end{proof}

\subsection{Explicit forms of the coefficients}

\noindent
Formula (\ref{coefib}) is not appropriate for numerical computations.
A more practical formula to compute the coefficients of the above
multi-point Taylor expansion is given in the following proposition.
First we have a definition:

\begin{definition}
\label{def2}
     Let $f(w)$ be analytic at $w$; then for $n=0,1,2,\ldots$
the differentential operator $D_w^n f(w)$ is defined by
\begin{equation}
D_w^n f(w)={1\over n!} {d^n\over dw^n} f(w).
\end{equation}
\end{definition}

\begin{prop}
\label{prop1}
     The coefficients $a_{n,j,l}$, for $n=1,2,3,\ldots$,
$j=1,2,\ldots,p$, $l=0,1,\ldots,m_j-1$
   in the expansion (\ref{expanib}) are also given by the formula:
\begin{equation}
     \label{coefibis}
\begin{split}
a_{n,j,l}=&
D_w^{nm_j+l}\left.\left[
{f(w)\over\prod_{s=1,s\ne j}^p(w-z_s)^{nm_s}}
\right]\right\vert_{w=z_j}+ \\ &
\sum_{k=1,k\ne j}^pD_w^{nm_k-1}
\left.\left[
{f(w)\over(w-z_j)^{l+1}\prod_{s=1,s\ne k}^p(w-z_s)^{nm_s}}
\right]\right\vert_{w=z_k}.
\end{split}
\end{equation}
\end{prop}

\begin{proof}
     We deform the contour of integration $C$ in equation
(\ref{coefib}) to any contour of the form
$C_1\cup C_2\cup\cdots\cup C_p$, also contained in $\Omega$,
where $C_k$, $k=1,2,\ldots,p$, is a simple closed loop which
encircles the
point $z_k$ in the counterclockwise direction and does not contain any
other point $z_j$, $j=1,2,\ldots,p$, $j\ne k$ inside (see Figure 3 (a)). Then,
\begin{equation}
\begin{split}
a_{n,j,l}= & {1\over 2\pi i}\sum_{k=1,k\ne j}^p
\int_{C_k}{f(w)\over
(w-z_j)^{l+1}\prod_{s=1,s\ne k}^p(w-z_s)^{nm_s}}
{dw\over (w-z_k)^{nm_k}}+  \\ &
{1\over 2\pi i}\int_{C_j}{f(w)\over
\prod_{s=1,s\ne j}^p(w-z_s)^{nm_s}}{dw\over (w-z_j)^{nm_j+l+1}},
\end{split}
\end{equation}
from which equation (\ref{coefibis}) follows.
\end{proof}

\begin{figure}[tb]
     \begin{center}\includegraphics[width=12cm]{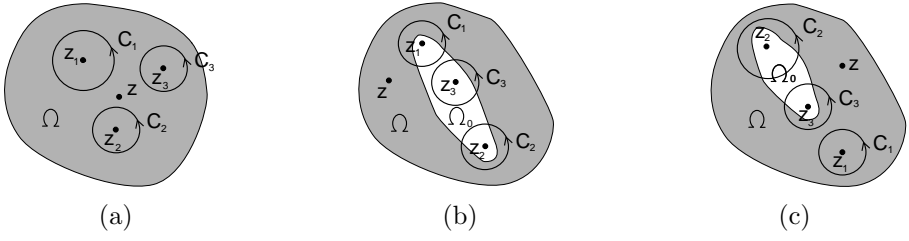}\end{center}
\centerline{(a) \hskip 4cm (b)\hskip 4cm (c)}
\caption{Integration contours $C_k$ for $p=3$
and $q=1$. (a) The function
$\prod_{s=1,s\ne k}^p(w-z_s)^{-nm_s}f(w)$ is analytic
inside $C_k$ for $k=1,2,\ldots,p$.
(b) The functions $\prod_{s=1,s\ne k}^p(w-z_s)^{-nm_s}g_k(w)$ and
$\prod_{s=1,s\ne k}^p(w-z_s)^{(n+1)m_s}g_k(w)$ are analytic inside
$C_k$ for $k=1,2,\ldots,p$.
(c) The functions $\prod_{s=1,s\ne k}^p(w-z_s)^{-nm_s}g_k(w)$,
$\prod_{s=1,s\ne k}^q(w-z_s)^{-nm_s}
\prod_{s=q+1,s\ne k}^p(w-z_s)^{nm_s}g_k(w)$ and
$\prod_{s=1,s\ne k}^q(w-z_s)^{-(n+1)m_s}\ \times\ $
$\prod_{s=q+1,s\ne k}^p(w-z_s)^{(n+1)m_s}g_k(w)$ are analytic inside
$C_k$ for $k=1,2,\ldots,p$.}
\label{Figure 3}
\end{figure}

\subsection{Multi-point Taylor polynomials}

\noindent
In Theorem \ref{theo1} we have assumed that the function $f(z)$ is analytic in
$\Omega$. If $f(z)$ is not analytic in
$\Omega$ but has a finite number of derivatives at $z_1$,
$z_2$,...,$z_p$, we can still define the multi-point Taylor polynomial of the
function $f(z)$ at $z_1$, $z_2$,...,$z_p$ in the following way:

\begin{definition}
\label{def3}
     Let $z$ be a real or complex variable. If $f(z)$ is
$Nm_k-1-$times differentiable at $z_1$, $z_2$,...,$z_p$,
we define the multi-point
Taylor polynomial of $f(z)$ at the points of $S$ and degree $mN-1$ as
\begin{equation}
P_N(z)\equiv\sum_{n=0}^{N-1}q_{n,m}(z)\prod_{k=1}^p(z-z_k)^{nm_k},
\end{equation}
where $q_{n,m}(z)$ is the polynomial of degree $m-1$
\begin{equation}
q_{n,m}(z)\equiv\sum_{j=1}^p
{\prod_{k=1,k\ne j}^p(z-z_k)^{m_k}\over \prod_{k=1,k\ne j}^p(z_j-z_k)^{m_k}}
\sum_{l=0}^{m_j-1}a_{n,j,l}(z-z_j)^l
\end{equation}
and the coefficients $a_{n,j,l}$ are given in (\ref{coefibis}).
\end{definition}

\begin{prop}
\label{prop2}
    In the conditions of the above definition, define
the remainder of the
approximation of $f(z)$ by $P_N(z)$ at the points of $S$ as
\begin{equation}
r_N(z)\equiv f(z)-P_N(z).
\end{equation}
Then, (i) $r_N(z)=o(z-z_k)^{Nm_k-1}$ as $z\to z_k$, $k=1,2,\ldots,p$.
(ii) If $f(z)$ is $Nm_k-$times differentiable at $z_k$ for some $k$, then
$r_N(z)=O(z-z_k)^{Nm_k}$ as $z\to z_k$.
\end{prop}

\begin{proof}
     The proof is trivial if $f(z)$ is analytic at every
$z_1$, $z_2$,...,$z_p$ by using (\ref{remib}).
In any case, for real or complex variable, the proof follows
by using l'H\^opital's rule and (\ref{coefibis}).
\end{proof}

\begin{remark}
\label{rem1}
     Observe that the Taylor polynomial of $f(z)$ at the points of $S$
and degree $mN-1$ is the Hermite's interpolation polynomial of
$f(z)$ at $z_1$, $z_2$,...,$z_p$ with data $f(z_k)$,
$f'(z_k)$,...,$f^{(Nm_k-1)}(z_k)$, $k=1,2,\ldots,p$.
\end{remark}

\section{Multi-point Laurent expansions}

\noindent
In the standard theory for Taylor and Laurent expansions much analogy
exists between the two types of expansions. For multi-point expansions, we
have a similar resemblance in the representations of coefficients and remainders.

\begin{theorem}
\label{theo3}
     Let $\Omega_0$ and $\Omega$ be closed and open sets,
respectively, of the complex plane, and $\Omega_0\subset\Omega\subset\Cs$. Let
$f(z)$ be an analytic function on $\Omega\setminus\Omega_0$
and $z_1$, $z_2$,...,$z_p\in\Omega_0$ (That is, $S\in\Omega_0$).
Then, for any $z\in\Omega\setminus\Omega_0$,
$f(z)$ has the multi-point Laurent expansion
\begin{equation}
     \label{expanii}
f(z)= \sum_{n=0}^{N-1}q_{n,m}(z)\prod_{k=1}^p(z-z_k)^{nm_k}+
\sum_{n=0}^{N-1}t_{n,m}(z)\prod_{k=1}^p(z-z_k)^{-(n+1)m_k}+ r_N(z),
\end{equation}
where $q_{n,m}(z)$ is the polynomial of degree $m-1$
\begin{equation}
     \label{polybis}
q_{n,m}(z)\equiv\sum_{j=1}^p
{\prod_{k=1,k\ne j}^p(z-z_k)^{m_k}\over \prod_{k=1,k\ne j}^p(z_j-z_k)^{m_k}}
\sum_{l=0}^{m_j-1}a_{n,j,l}(z-z_j)^l
\end{equation}
and the coefficients $a_{n,j,l}$ of this polynomial are given by the
Cauchy integral
\begin{equation}
     \label{coefii}
a_{n,j,l}\equiv{1\over 2\pi i}
\int_{\Gamma_1}{f(w)\,dw\over (w-z_j)^{l+1}\prod_{k=1}^p(w-z_k)^{nm_k}}.
\end{equation}
Also, $t_{n,m}(z)$ is the polynomial of degree $m-1$
\begin{equation}
     \label{polyi}
t_{n,m}(z)\equiv\sum_{j=1}^p
{\prod_{k=1,k\ne j}^p(z-z_k)^{m_k}\over \prod_{k=1,k\ne j}^p(z_j-z_k)^{m_k}}
\sum_{l=0}^{m_j-1}b_{n,j,l}(z-z_j)^l,
\end{equation}
where the coefficients $b_{n,j,l}$ of this polynomial are given by the
Cauchy integral
\begin{equation}
     \label{coefiii}
b_{n,j,l}\equiv{1\over 2\pi i}
\int_{\Gamma_2}\prod_{k=1}^p(w-z_k)^{m_k(n+1)}{f(w)\,dw\over (w-z_j)^{l+1}}.
\end{equation}
The remainder term $r_N(z)$ is given by the Cauchy integrals
\begin{equation}
     \label{remii}
\begin{split}
r_N(z)\equiv & {1\over 2\pi i}\int_{\Gamma_1}{f(w)\,dw\over (w-z)
\prod_{k=1}^p(w-z_k)^{Nm_k}}\prod_{k=1}^p(z-z_k)^{Nm_k}- \\ &
{1\over 2\pi i}\int_{\Gamma_2}\prod_{k=1}^p(w-z_k)^{Nm_k}{f(w)\,dw\over w-z
}{1\over\prod_{k=1}^p(z-z_k)^{Nm_k}}.
\end{split}
\end{equation}
In these integrals, the contours of integration ${\Gamma}_1$ and ${\Gamma}_2$
are simple closed loops contained in $\Omega\setminus\Omega_0$ which
encircle the
points $z_1$, $z_2$,...,$z_p$ in the counterclockwise direction.
Moreover, $z$ is not inside ${\Gamma}_2$,
whereas ${\Gamma}_1$ encircles ${\Gamma}_2$ and the
point $z$ (see Figure 4 (a)).

\noindent
The expansion (\ref{expanii}) is convergent for $z$ inside the
"lemniscate annulus" (see Figure 5)
\begin{equation}
     \label{domainii}
A_{p}\equiv\lbrace z\in\Omega\setminus\Omega_0, \hskip 2mm
r_2<\prod_{k=1}^p\vert z-z_k\vert^{m_k}<r_1\rbrace,
\end{equation}
where
\begin{equation}
     \label{erei}
r_1\equiv {\rm Inf}_{w\in\Css\setminus\Omega}\left\lbrace
\prod_{k=1}^p\vert w-z_k\vert^{m_k}\right\rbrace,\quad
r_2\equiv{\rm Sup}_{w\in\Omega_0}\left\lbrace
\prod_{k=1}^p\vert w-z_k\vert^{m_k}\right\rbrace.
\end{equation}
\end{theorem}

\begin{proof}
     By Cauchy's theorem,
\begin{equation}
     \label{cauchybis}
f(z)={1\over 2\pi i}\int_{{\Gamma}_1}{f(w)dw\over w-z}-
{1\over 2\pi i}\int_{{\Gamma}_2}{f(w)\,dw\over w-z},
\end{equation}
where ${\Gamma}_1$ and ${\Gamma}_2$ are the contours defined above.
First we assume that the $m$ points of the set
$S$ are all distinct and later we will let the first $m_1$
points to coalesce to $z_1$, the second $m_2$ points to coalesce to
$z_2$ and so on. We
substitute (\ref{lequal})-(\ref{ui}) into the first integral above and
\begin{equation}
{1\over w-z}=-{H_m(w,z;z_1,\ldots,z_m)\over \prod_{k=1}^m(z-z_k)}{1\over 1-u},
\hskip 1cm
u\equiv {\prod_{k=1}^m(w-z_k)\over \prod_{k=1}^m(z-z_k)},
\end{equation}
where $H_m(w,z;z_1,\ldots,z_m)$ is defined in (\ref{defache}),
into the second one. Now we introduce the expansion (\ref{expanu}) of 
the factor
$(1-u)^{-1}$ in both integrals in (\ref{cauchybis}).
Using (\ref{ache}) and after straightforward calculations we obtain
\begin{equation}
f(z)= \sum_{n=0}^{N-1}q_{n,m}(z)\prod_{k=1}^m(z-z_k)^{n}+
\sum_{n=0}^{N-1}t_{n,m}(z)\prod_{k=1}^m(z-z_k)^{-n-1}+ r_N(z),
\end{equation}
where $q_{n,m}(z)$ is given by formulas (\ref{poly}) and (\ref{coefi})
replacing the contour $C$ by $\Gamma_1$. Also,
\begin{equation}
t_{n,m}(z)=\sum_{j=1}^m b_{n,j}{\prod_{k=1,k\ne j}^m(z-z_k)
\over \prod_{k=1,k\ne j}^m(z_j-z_k)},
\hskip 1cm
b_{n,j}\equiv{1\over 2\pi i}
\int_{\Gamma_2}\prod_{k=1}^m(w-z_k)^{n+1}{f(w)\,dw\over w-z_j},
\end{equation}
and
\begin{equation}
\begin{split}
r_N(z)= & {1\over 2\pi i}\int_{\Gamma_1}{f(w)\,dw\over (w-z)
\prod_{k=1}^m(w-z_k)^N}\prod_{k=1}^m(z-z_k)^N- \\ &
{1\over 2\pi i}\int_{\Gamma_2}\prod_{k=1}^m(w-z_k)^N{f(w)\,dw\over w-z
}{1\over\prod_{k=1}^m(z-z_k)^N}.
\end{split}
\end{equation}
Now we write
\begin{equation}
t_{n,m}(z)={1\over 2\pi i}
\int_{\Gamma_2}\prod_{k=1}^m(w-z_k)^{n+1}f(w)\,dw\sum_{j=1}^m
{\prod_{k=1,k\ne j}^m(z-z_k)\over (w-z_j)\prod_{k=1,k\ne j}^m(z_j-z_k)}
\end{equation}
and repeat the steps following (\ref{cuco}) in the proof of Theorem \ref{theo2}
for $q_{n,m}(z)$ and $t_{n,m}(z)$.

For any $z\in A_p$ we can take simple closed loops ${\Gamma}_1$ and
${\Gamma}_2$ in $\Omega\setminus\Omega_0$ such that (see Figure 4 (b))
\begin{equation}
\prod_{k=1}^p\vert z-z_k\vert^{m_k}<
\prod_{k=1}^p\vert w-z_k\vert^{m_k} \hskip 1cm \forall \hskip 2mm
w\in{\Gamma}_1
\end{equation}
and
\begin{equation}
\prod_{k=1}^p\vert z-z_k\vert^{m_k}>
\prod_{k=1}^p\vert w-z_k\vert^{m_k} \hskip 1cm
\forall \hskip 2mm w\in{\Gamma}_2.
\end{equation}
On these contours $\vert f(w)\vert$ is bounded by some constant $C$: $\vert
f(w)\vert\le C$. Introducing these bounds in (\ref{remii}) we see that
$\lim_{N\to\infty}r_N(z)=0$ and the proof follows.
\end{proof}

\begin{figure}[tb]
     \begin{center}\includegraphics[width=10cm]{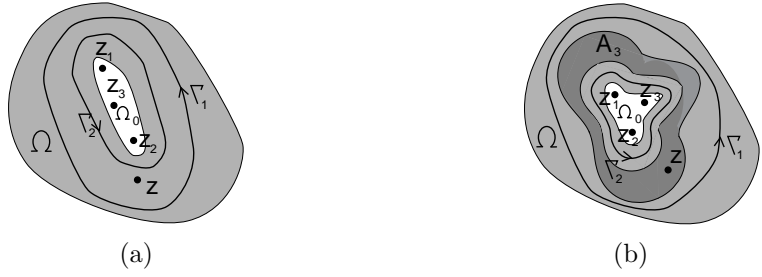}\end{center}
\centerline{(a) \hskip 6cm (b)}
\caption{The case $p=3$.
(a) Contours $\Gamma_1$ and $\Gamma_2$ in the integrals
(\ref{coefii}), (\ref{coefiii}) and (\ref{remii}).
(b) For $z\in A_{p}$, we can take a contour $\Gamma_2$ in
$\Omega$ located between $\Omega_0$ and $A_{p}$ and a contour
$\Gamma_1$ in $\Omega$ such that $A_{p}$ is inside this contour. Therefore,
$\prod_{k=1}^p\vert z-z_k\vert^{m_k}<\prod_{k=1}^p\vert w-z_k\vert^{m_k}$
$\forall$ $w\in\Gamma_1$ and
$\prod_{k=1}^p\vert w-z_k\vert^{m_k}<\prod_{k=1}^p\vert z-z_k\vert^{m_k}$
$\forall$ $w\in\Gamma_2$.}
\label{Figure 4}
\end{figure}

\begin{figure}[tb]
     \begin{center}\includegraphics[width=10cm]{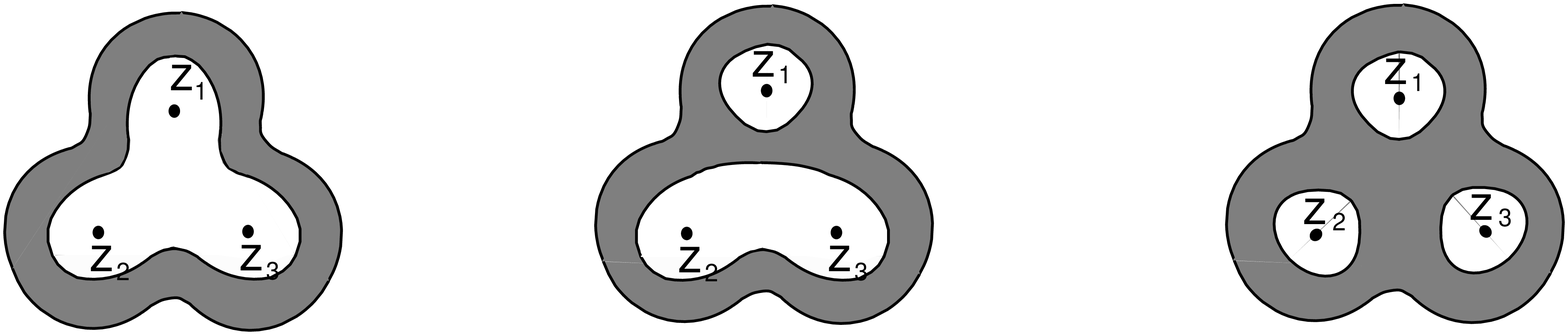}\end{center}
\centerline{$(r_1^a,r_2^a)$ \hskip 3cm
$(r_1^a,r_2^b)$ \hskip 3cm $(r_1^a,r_2^c)$}
\bigskip
     \begin{center}\includegraphics[width=10cm]{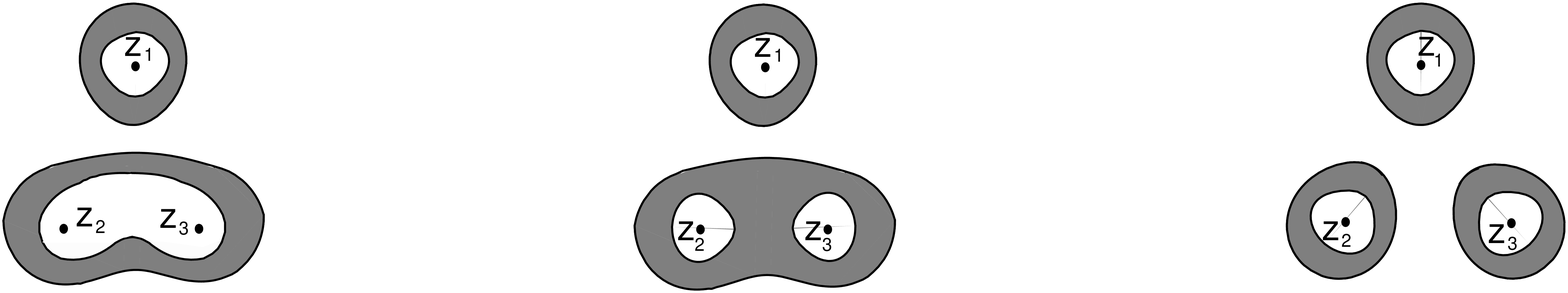}\end{center}
\centerline{$(r_1^b,r_2^b)$ \hskip 3cm
$(r_1^b,r_2^c)$ \hskip 3cm $(r_1^c,r_2^c)$}
\caption{Shape of the "lemniscate annulus" $A_p$ for $p=3$.
It depends on the relative size of the parameters $r_1$ and $r_2$
defined in (\ref{erei}). The different forms are labeled by $(r_1,r_2)$
with $r_1>r_2$. In these pictures
$\vert z_2-z_3\vert<\vert z_1-z_3\vert$, $\vert z_1-z_2\vert$
   and $r_2^a>r_2^b>r_2^c$.}
\label{Figure 5}
\end{figure}

If the only singularities of $f(z)$ inside $\Omega_0$ are just poles at
$z_1$, $z_2$,...,$z_p$, then
alternative formulas of (\ref{coefii}) and (\ref{coefiii}) for computing the
coefficients of the above multi-point
Laurent expansion is given in the following proposition.

\begin{prop}
\label{prop3}
     Suppose that $g_k(z)\equiv (z-z_k)^{\rho_k}f(z)$,
$k=1,2,\ldots,p$ are analytic functions in $\Omega$ for certain
$\rho_1$, $\rho_2$,...,$\rho_k\in\Ns$. Then, for $n=0,1,2,\ldots$,
coefficients $a_{n,j,l}$ and $b_{n,j,l}$ in expansion (\ref{expanii}) are
also given by the formulas:
\begin{equation}
     \label{coefiibis}
\begin{split}
a_{n,j,l}= &
\left.\sum_{k=1,k\ne j}^p
D_w^{nm_k+\rho_k-1}\left[
{g_k(w)\over(w-z_j)^{l+1}\prod_{s=1,s\ne k}^p(w-z_s)^{nm_s}}\right]
\right\vert_{w=z_k} + \\ &
\quad D_w^{nm_j+\rho_j+l}
\left. \left[
{g_j(w)\over\prod_{s=1,s\ne j}^p(w-z_s)^{nm_s}}
\right]\right\vert_{w=z_j}
\end{split}
\end{equation}
and
\begin{equation}
     \label{coefiibisbis}
\begin{split}
b_{n,j,l}= &
\left.\sum_{k=1,k\ne j}^p
D_w^{\rho_k-(n+1)m_k-1}
\left[
{g_k(w)\over(w-z_j)^{l+1}}\prod_{s=1,s\ne k}^p(w-z_s)^{(n+1)m_s}
\right]\right\vert_{w=z_k}+ \\ &
\quad
D_w^{\rho_j-(n+1)m_j+l}\left. \left[
{g_j(w)\prod_{s=1,s\ne j}^p(w-z_s)^{(n+1)m_s}}
\right]\right\vert_{w=z_j}.
\end{split}
\end{equation}
\end{prop}

\begin{proof}
     We deform both contours ${\Gamma}_1$ and ${\Gamma}_2$
of equations (\ref{coefii}) and (\ref{coefiii}), respectively, to any
contour of the form $C_1\cup C_2\cup\cdots\cup C_p$
contained in $\Omega$, where
$C_k$, $k=1,2,\ldots,p$ is a simple closed loop which
encircles the point $z_k$ in the counterclockwise
direction and does not contain the point $z_j$ $j=1,2,\ldots,p$, $j\ne k$
inside (see Figure 3 (b)). Then,
\begin{equation}
\begin{split}
a_{n,j,l}= & {1\over 2\pi i}
\sum_{k=1,k\ne j}^p\int_{C_k}{g_k(w)\over
(w-z_j)^{l+1}\prod_{s=1,s\ne k}^p(w-z_s)^{nm_s}}
{dw\over (w-z_k)^{nm_k+\rho_k}}+ \\ &
{1\over 2\pi i}\int_{C_j}{g_j(w)\over
\prod_{s=1,s\ne j}^p(w-z_s)^{nm_s}}
{dw\over (w-z_j)^{nm_j+\rho_j+l+1}}
\end{split}
\end{equation}
and
\begin{equation}
\begin{split}
b_{n,j,l}= & {1\over 2\pi i}
\sum_{k=1,k\ne j}^p\int_{C_k}
{\prod_{s=1,s\ne k}^p(w-z_s)^{(n+1)m_s}\over
(w-z_j)^{l+1}}{g_k(w)dw\over (w-z_k)^{\rho_k-(n+1)m_k}}+  \\ &
{1\over 2\pi i}\int_{C_j}\prod_{s=1,s\ne j}^p(w-z_s)^{(n+1)m_s}
{g_j(w)dw\over (w-z_j)^{\rho_j-(n+1)m_j+l+1}}.
\end{split}
\end{equation}

  From here, equations (\ref{coefiibis}) and (\ref{coefiibisbis}) follow.
\end{proof}

\begin{remark}
\label{rem2}
     Let $z$ be a real or complex variable.
Suppose that $g_k(z)\equiv (z-z_k)^{\rho_k}f(z)$
is $\rho_k-1-$times differentiable at every $z_k$ in $S$
for some $\rho_k\in\Ns$, $k=1,2,\ldots,p$. Define
\begin{equation}
g(z)\equiv f(z)-
\sum_{n=0}^{M}t_{n,m}(z)\prod_{k=1}^p(z-z_k)^{-(n+1)m_k},
\end{equation}
where $M\equiv\lfloor$Max$\lbrace (\rho_1-1)/m_1,(\rho_2-1)/m_2,\ldots,
(\rho_p-1)/m_p\rbrace\rfloor$ and $t_{n,m}(z)$
is the polynomial defined in (\ref{polyi}) and (\ref{coefiibisbis}).
Then, the thesis of Proposition \ref{prop2}
holds for $f(z)$ replaced by $g(z)$. Moreover, if
$\prod_{k=1}^p(z-z_k)^{\rho_k}f(z)$ is an
analytic function in $\Omega$, then the thesis of Theorem \ref{theo2}
applies to $g(z)$.
\end{remark}

\section{Multi-point Taylor-Laurent expansions}

\noindent
For multi-point expansions we have the possibility (that we do not
have in the standard theory) of expanding in Taylor series in some points and
in Laurent series in other points.

\begin{theorem}
\label{theo4}
     Let $\Omega_0$ and $\Omega$ be closed and open
sets, respectively, of the complex
plane, and $\Omega_0\subset\Omega\subset\Cs$. Let $f(z)$
be an analytic function on $\Omega\setminus\Omega_0$,
$z_1,z_2,\ldots,z_q\in\Omega\setminus\Omega_0$ and
$z_{q+1},z_{q+2},\ldots,z_p\in\Omega_0$ ($q$ points are in
$\Omega\setminus\Omega_0$ and $p-q$ points are in $\Omega_0$). Write
$s\equiv m_1+m_2+\cdots+m_q$. Then, for
$z\in\Omega\setminus\Omega_0$, $f(z)$ has the Taylor-Laurent expansion
\begin{equation}
     \label{expaniii}
\begin{split}
f(z)= & \sum_{n=0}^{N-1}q_{n,m}(z)\prod_{k=1}^p(z-z_k)^{nm_k}+
\sum_{n=0}^{N-1}t_{n,m}^{(1)}(z){\prod_{k=1}^q(z-z_k)^{nm_k}\over
\prod_{k=q+1}^p(z-z_k)^{nm_k}}+  \\ &
\sum_{n=0}^{N-1}t_{n,m}^{(2)}(z){\prod_{k=1}^q(z-z_k)^{(n+1)m_k}\over
\prod_{k=q+1}^p(z-z_k)^{(n+1)m_k}}+ r_N(z),
\end{split}
\end{equation}
where $q_{n,m}(z)$ is the polynomial of degree $m-1$
\begin{equation}
     \label{polybisbis}
q_{n,m}(z)\equiv\sum_{j=1}^p
{\prod_{k=1,k\ne j}^p(z-z_k)^{m_k}\over \prod_{k=1,k\ne j}^p(z_j-z_k)^{m_k}}
\sum_{l=0}^{m_j-1}a_{n,j,l}(z-z_j)^l
\end{equation}
and the coefficients $a_{n,j,l}$ of this polynomial are given by the
Cauchy integral
\begin{equation}
     \label{coefiv}
a_{n,j,l}\equiv{1\over 2\pi i}
\int_{\Gamma_1}{f(w)\,dw\over (w-z_j)^{l+1}\prod_{k=1}^p(w-z_k)^{nm_k}}.
\end{equation}
Also, $t_{n,m}^{(1)}(z)$ and $t_{n,m}^{(2)}(z)$
are the following polynomials of degrees $s-1$ and $m-s-1$ respectivelly,
\begin{equation}
     \label{polyii}
t_{n,m}^{(1)}(z)\equiv -\sum_{j=1}^q
{\prod_{k=1,k\ne j}^q(z-z_k)^{m_k}\over \prod_{k=1,k\ne j}^q(z_j-z_k)^{m_k}}
\sum_{l=0}^{m_j-1}b_{n,j,l}(z-z_j)^l,
\end{equation}
where the coefficients $b_{n,j,l}$ of this polynomial are given by the
Cauchy integral
\begin{equation}
     \label{coefv}
b_{n,j,l}\equiv{1\over 2\pi i}
\int_{\Gamma_2}{\prod_{k=q+1}^p(w-z_k)^{nm_k}\over
\prod_{k=1}^q(w-z_k)^{nm_k}}
{f(w)\,dw\over (w-z_j)^{l+1}}.
\end{equation}
\begin{equation}
     \label{polyiii}
t_{n,m}^{(2)}(z)\equiv \sum_{j=q+1}^p
{\prod_{k=q+1,k\ne j}^p(z-z_k)^{m_k}\over
\prod_{k=q+1,k\ne j}^p(z_j-z_k)^{m_k}}
\sum_{l=0}^{m_j-1}c_{n,j,l}(z-z_j)^l,
\end{equation}
where the coefficients $c_{n,j,l}$ of this polynomial are given by the
Cauchy integral
\begin{equation}
     \label{coefvi}
c_{n,j,l}\equiv{1\over 2\pi i}
\int_{\Gamma_2}{\prod_{k=q+1}^p(w-z_k)^{(n+1)m_k}\over
\prod_{k=1}^q(w-z_k)^{(n+1)m_k}}
{f(w)\,dw\over (w-z_j)^{l+1}}.
\end{equation}
The remainder term $r_N(z)$ is given by the Cauchy integrals
\begin{equation}
     \label{remiii}
\begin{split}
r_N(z)\equiv & {1\over 2\pi i}\int_{\Gamma_1}{f(w)\,dw\over (w-z)
\prod_{k=1}^p(w-z_k)^{Nm_k}}\prod_{k=1}^p(z-z_k)^{Nm_k}- \\ &
{1\over 2\pi i}\int_{\Gamma_2}{\prod_{k=q+1}^p(w-z_k)^{Nm_k}\over
\prod_{k=1}^q(w-z_k)^{Nm_k}}{f(w)dw\over w-z
}{\prod_{k=1}^q(z-z_k)^{Nm_k}\over\prod_{k=q+1}^p(z-z_k)^{Nm_k}}.
\end{split}
\end{equation}
In these integrals, the contours of integration ${\Gamma}_1$ and ${\Gamma}_2$
are simple closed loops
contained in $\Omega\setminus\Omega_0$ which encircle $\Omega_0$ in the
counterclockwise direction. Moreover,
the points $z$ and $z_1$, $z_2$,...,$z_q$ are not inside ${\Gamma}_2$,
whereas ${\Gamma}_1$ encircles ${\Gamma}_2$ and the points
$z$ and $z_1$, $z_2$,...,$z_q$ (see Figure 6 (a)).

\noindent
The expansion (\ref{expaniii}) is convergent in the region (Figure 7)
\begin{equation}
     \label{domainiii}
\begin{split}
D_{q,p}\equiv\biggl\lbrace z\in\Omega\setminus\Omega_0, & \hskip 1mm
\prod_{k=1}^p\vert(z-z_k)\vert^{m_k}<r_1, \\ & \left.
\prod_{k=1}^q\vert(z-z_k)\vert^{m_k}<
r_2\prod_{k=q+1}^p\vert(z-z_k)\vert^{m_k}\right\rbrace
\end{split}
\end{equation}
where $r_1\equiv$ Inf$_{w\in\Css\setminus\Omega}\left\lbrace
\prod_{k=1}^p\vert(w-z_k)\vert^{m_k}\right\rbrace$
and

\noindent
$r_2\equiv$ Inf$_{w\in\Omega_0}\left\lbrace
\prod_{k=1}^q\vert(w-z_k)\vert^{m_k}
\prod_{k=q+1}^p\vert(w-z_k)^{-1}\vert^{m_k}\right\rbrace$.
\end{theorem}

\begin{proof}
     By Cauchy's theorem,
\begin{equation}
     \label{cauchybisbis}
f(z)={1\over 2\pi i}\int_{{\Gamma}_1}{f(w)\,dw\over w-z}-
{1\over 2\pi i}\int_{{\Gamma}_2}{f(w)\,dw\over w-z},
\end{equation}
where ${\Gamma}_1$ and ${\Gamma}_2$ are the contours defined above.

First we assume that the $m$ points of the set
$S$ are all distinct. Later we will let the first $m_1$
points coalesce to $z_1$, the second $m_2$ points to
$z_2$, and so on. We
substitute (\ref{lequal})-(\ref{ui}) into the first integral above and
\begin{equation}
{1\over w-z}={F_m(w,z;z_1,\ldots,z_m)\over
\prod_{k=1}^s(w-z_k)\prod_{k=s+1}^m(z-z_k)}{1\over 1-u},
\end{equation}
where
\begin{equation}
u\equiv {\prod_{k=1}^s(z-z_k)\prod_{k=s+1}^m(w-z_k)\over
\prod_{k=1}^s(w-z_k)\prod_{k=s+1}^m(z-z_k)}
\end{equation}
and
\begin{equation}
\begin{split}
F_m(w,z;z_1,\ldots,z_m)\equiv {1\over w-z}&\left[
\prod_{k=1}^s(w-z_k)\prod_{k=s+1}^m(z-z_k)- \right. \\ & \left.
\prod_{k=1}^s(z-z_k)\prod_{k=s+1}^m(w-z_k)\right]
\end{split}
\end{equation}
into the second one. Next we introduce the expansion (\ref{expanu}) 
of the factor
$(1-u)^{-1}$ in both integrals in (\ref{cauchybisbis}). We observe that
$F_m(w,z;z_1,\ldots,z_m)$ may be written as
\begin{equation}
\begin{split}
F_m(w,z;z_1,\ldots,z_m)=
&H_s(w,z;z_1,\ldots,z_p)\prod_{k=s+1}^m(z-z_k) -\\
&H_{m-s}(w,z;z_{p+1},\ldots,z_m)\prod_{k=1}^s(z-z_k),
\end{split}
\end{equation}
where $H_m(w,z;z_1,\ldots,z_m)$ is defined in (\ref{defache}).
Using this decomposition,
equation (\ref{ache}), and after straightforward calculations we obtain
\begin{equation}
\begin{split}
f(z)= & \sum_{n=0}^{N-1}q_{n,m}(z)\prod_{k=1}^m(z-z_k)^{n}+
\sum_{n=0}^{N-1}t_{n,m}^{(1)}(z){\prod_{k=1}^s(z-z_k)^{n}
\over\prod_{k=s+1}^m(z-z_k)^{n}}+ \\ &
\sum_{n=0}^{N-1}t_{n,m}^{(2)}(z){\prod_{k=1}^s(z-z_k)^{n+1}
\over\prod_{k=s+1}^m(z-z_k)^{n+1}}+ r_N(z),
\end{split}
\end{equation}
where $q_{n,m}(z)$ is given by formulas (\ref{poly}) and (\ref{coefi})
replacing the contour $C$ by $\Gamma_1$. Also,
\begin{equation}
t_{n,m}^{(1)}(z)=-\sum_{j=1}^s b_{n,j}{\prod_{k=1,k\ne j}^s(z-z_k)
\over \prod_{k=1,k\ne j}^s(z_j-z_k)}
\end{equation}
with
\begin{equation}
b_{n,j}\equiv{1\over 2\pi i}
\int_{\Gamma_2}{\prod_{k=s+1}^m(w-z_k)^n\over
\prod_{k=1}^s(w-z_k)^n}{f(w)\,dw\over w-z_j},
\end{equation}
\begin{equation}
t_{n,m}^{(2)}(z)=\sum_{j=s+1}^m c_{n,j}{\prod_{k=s+1,k\ne j}^m(z-z_k)
\over \prod_{k=s+1,k\ne j}^m(z_j-z_k)}
\end{equation}
with
\begin{equation}
c_{n,j}\equiv{1\over 2\pi i}
\int_{\Gamma_2}{\prod_{k=s+1}^m(w-z_k)^{n+1}\over
\prod_{k=1}^s(w-z_k)^{n+1}}{f(w)\,dw\over w-z_j}
\end{equation}
and
\begin{equation}
\begin{split}
r_N(z)= & {1\over 2\pi i}\int_{\Gamma_1}{f(w)\,dw\over (w-z)
\prod_{k=1}^m(w-z_k)^{N}}\prod_{k=1}^m(z-z_k)^{N}- \\ &
{1\over 2\pi i}\int_{\Gamma_2}{\prod_{k=s+1}^m(w-z_k)^{N}
\over\prod_{k=1}^s(w-z_k)^{N}}{f(w)\,dw\over w-z
}{\prod_{k=1}^s(z-z_k)^{N}\over\prod_{k=s+1}^m(z-z_k)^{N}}.
\end{split}
\end{equation}
Now we write
\begin{equation}
\begin{split}
t_{n,m}^{(1)}(z)
&=-{1\over 2\pi i}
\int_{\Gamma_2}{\prod_{k=s+1}^m(w-z_k)^n
\over\prod_{k=1}^s(w-z_k)^n}
\sum_{j=1}^s
{f(w)\prod_{k=1,k\ne j}^s(z-z_k)\over (w-z_j)\prod_{k=1,k\ne j}^s(z_j-z_k)}
\,dw,\\
t_{n,m}^{(2)}(z)
&={1\over 2\pi i}
\int_{\Gamma_2}{\prod_{k=s+1}^m(w-z_k)^{n+1}
\over\prod_{k=1}^s(w-z_k)^{n+1}}
   \sum_{j=s+1}^m
{f(w)\prod_{k=s+1,k\ne j}^m(z-z_k)\over (w-z_j)\prod_{k=s+1,k\ne j}^m(z_j-z_k)}
\,dw,
\end{split}
\end{equation}
and repeat the steps following (\ref{cuco}) in Theorem \ref{theo2}
for $q_{n,m}(z)$,
$t_{n,m}^{(1)}(z)$ and $t_{n,m}^{(2)}(z)$.

For any $z$ verifying (\ref{domainiii}), we can take simple closed loops
${\Gamma}_1$ and ${\Gamma}_2$ in $\Omega\setminus\Omega_0$ such that
  (see Figure 6 (b))
\begin{equation}
\prod_{k=1}^p\vert z-z_k\vert^{m_k}<
\prod_{k=1}^p\vert w-z_k\vert^{m_k} \hskip 1cm \forall \hskip 2mm
w\in{\Gamma}_1
\end{equation}
and
\begin{equation}
\prod_{k=1}^q\vert w-z_k\vert^{m_k}
\prod_{k=q+1}^p\vert z-z_k\vert^{m_k}>
\prod_{k=1}^q\vert z-z_k\vert^{m_k}
\prod_{k=q+1}^p\vert w-z_k\vert^{m_k}
\end{equation}
$\forall$ $w\in{\Gamma}_2$.
On these contours $\vert f(w)\vert$ is bounded by some constant $C$: $\vert
f(w)\vert\le C$. Introducing these bounds in (\ref{remiii}) we see that
$\lim_{N\to\infty}r_N(z)=0$ and the proof follows.
\end{proof}

\begin{figure}[tb]
     \begin{center}\includegraphics[width=10cm]{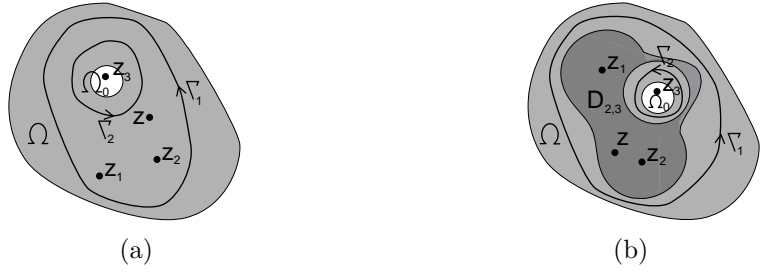}\end{center}
	\centerline{(a) \hskip 6cm (b)}
\caption{The case $q=2$, $p=3$.
(a) Contours $\Gamma_1$ and $\Gamma_2$ in the integrals
(\ref{coefiv}), (\ref{coefv}), (\ref{coefvi}) and (\ref{remiii}).
   (b) For $z\in D_{q,p}$, we can take a contour $\Gamma_2$
located between $\Omega_0$ and $D_{q,p}$ and a contour $\Gamma_1$
in $\Omega$ with $D_{q,p}$ inside this contour.
Therefore,
$\prod_{k=1}^p\vert z-z_k\vert^{m_k}<\prod_{k=1}^p\vert w-z_k\vert^{m_k}$
$\forall$ $w\in\Gamma_1$ and
$\prod_{k=1}^q\vert z-z_k\vert^{m_k}\prod_{k=q+1}^p\vert w-z_k\vert^{m_k}<
\prod_{k=1}^q\vert w-z_k\vert^{m_k}\prod_{k=q+1}^p\vert z-z_k\vert^{m_k}$
$\forall$ $w\in\Gamma_2$.}
\label{Figure 6}
\end{figure}

\begin{figure}[tb]
     \begin{center}\includegraphics[width=11cm]{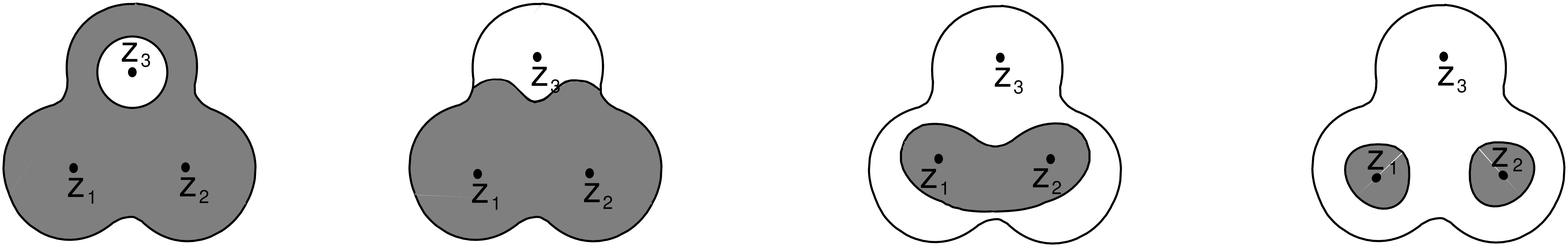}\end{center}
	\centerline{$(r_1^a,r_2^a)$ \hskip 1.8cm
$(r_1^a,r_2^b)$  \hskip 1.8cm $(r_1^a,r_2^c)$  \hskip 1.8cm
$(r_1^a,r_2^d)$}
\bigskip
     \begin{center}\includegraphics[width=11cm]{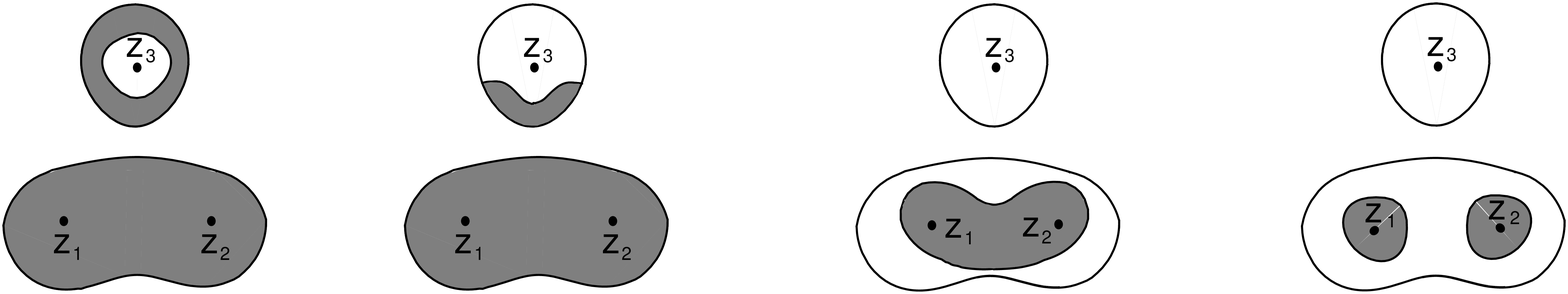}\end{center}
	\centerline{$(r_1^b,r_2^a)$  \hskip 1.8cm
$(r_1^b,r_2^b)$  \hskip 1.8cm $(r_1^b,r_2^c)$  \hskip 1.8cm
$(r_1^b,r_2^d)$}
\bigskip
     \begin{center}\includegraphics[width=11cm]{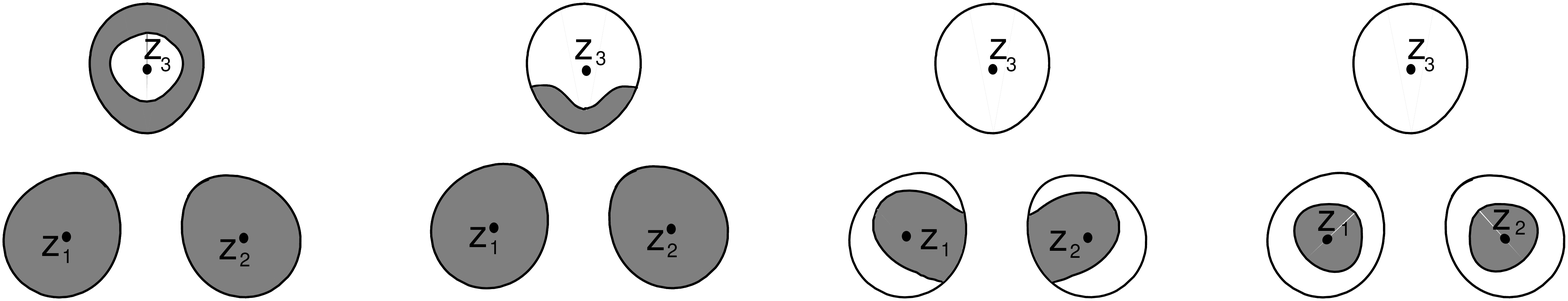}\end{center}
	\centerline{$(r_1^c,r_2^a)$ \hskip 1.8cm
$(r_1^c,r_2^b)$ \hskip 1.8cm $(r_1^c,r_2^c)$ \hskip 1.8cm
$(r_1^c,r_2^d)$}
\caption{The region $D_{q,p}$ defined in Theorem \ref{theo4} is given by
$D_{q,p}=O_p\bigcap B_{q,p}$, where $O_p$
is the "lemniscate domain" of foci $z_1$, $\ldots$, $z_p$ and parameter
$r_1$. Also,
$B_{q,p}\equiv \lbrace z\in\Cs$, $\prod_{k=1}^q\vert(z-z_k)\vert^{m_k}<
r_2\prod_{k=q+1}^p\vert(z-z_k)\vert^{m_k}\rbrace$. This pictures show the
topologically different forms of $D_{q,p}$ depending on the
relative value of $r_1$ and $r_2$ when $q=2$ and $p=3$. The pictures are
labeled with $(r_1,r_2)$. In these pictures $z_1$
$\vert z_1-z_2\vert<\vert z_1-z_3\vert$, $\vert z_2-z_3\vert$
   and $r_2^a>r_2^b>r_2^c>r_2^d$.}
\label{Figure 7}
\end{figure}

If the only singularities of $f(z)$ inside $\Omega_0$ are just poles at
$z_{q+1}$, $z_{q+2}$,...,$z_p$, then alternative formulas of
(\ref{coefiv}), (\ref{coefv}) and (\ref{coefvi}) for computing the coefficients
of the above two-point
Taylor-Laurent expansion is given in the following proposition.

\begin{prop}
\label{prop4}
     Suppose that $g_k(z)\equiv
(z-z_k)^{\rho_k}f(z)$ is an
analytic function in $\Omega$ for certain $\rho_k\in\Ns$ and
$k=q+1,q+2,\ldots,p$.
Define $g_k(w)=f(w)$ for $k=1,2,3,\ldots q$.
Then the coefficients $a_{n,j,l}$, $b_{n,j,l}$ and $c_{n,j,l}$ in the expansion
(\ref{expaniii}) are also given by the formulas:
\begin{equation}
     \label{coefvx}
\begin{split}
a_{n,j,l}= &
\left.\sum_{k=1,k\ne j}^q
D_w^{nm_k-1}\left[
{f(w)\over(w-z_j)^{l+1}\prod_{s=1,s\ne k}^p(w-z_s)^{nm_s}}
\right]\right\vert_{w=z_k}+ \\ &
\quad
\left.\sum_{k=q+1,k\ne j}^p
D_w^{nm_k+\rho_k-1}\left[
{g_k(w)\over(w-z_j)^{l+1}\prod_{s=1,s\ne k}^p(w-z_s)^{nm_s}}
\right]\right\vert_{w=z_k}+ \\ &
\quad
D_w^{nm_j+\rho_j+l}\left.
\left[
{g_j(w)\over\prod_{s=1,s\ne j}^p(w-z_s)^{nm_s}}\right]\right\vert_{w=z_j},
\end{split}
\end{equation}
\begin{equation}
     \label{coefiiibis}
b_{n,j,l}=
\left.\sum_{k=q+1}^p
D_w^{\rho_k-nm_k-1}\left[
{g_k(w)\prod_{s=q+1,s\ne k}^p(w-z_s)^{nm_s}\over(w-z_j)^{l+1}
\prod_{s=1}^q(w-z_s)^{nm_s}}\right]\right\vert_{w=z_k}.
\end{equation}
\begin{equation}
     \label{coefiiibisbis}
\begin{split}
c_{n,j,l}= &
\left.\sum_{k=q+1,k\ne j}^p
D_w^{\rho_k-(n+1)m_k-1}\left[
{g_k(w)\prod_{s=q+1,s\ne k}^p(w-z_s)^{(n+1)m_s}\over(w-z_j)^{l+1}
\prod_{s=1}^q(w-z_s)^{(n+1)m_s}}\right]
\right\vert_{w=z_k}+  \\ &
\quad
D_w^{\rho_j-(n+1)m_j+l}\left. \left[
{g_j(w)\prod_{s=q+1,s\ne j}^p(w-z_s)^{(n+1)m_s}
\over\prod_{s=1}^q(w-z_s)^{(n+1)m_s}}\right]
\right\vert_{w=z_j}.
\end{split}
\end{equation}
\end{prop}

\begin{proof}
     We deform both contours ${\Gamma}_1$ in equation (\ref{coefiv}) and
the  contour ${\Gamma}_2$ in equations (\ref{coefvx}) and (\ref{coefvi}) 
into any contour  of the form
${C_1}\cup{C_2}\cup\cdots\cup{C_p}$ contained in
$\Omega$, where ${C_k}$, $k=1,2,...,p$ is a simple closed loop
which encircles the point $z_k$ in the counterclockwise
direction with  $z_j$ not inside ${C_k}$,
$j=1,2,\ldots,p$, $j\ne k$ (see Figure 3 (c)). Then,
\begin{equation}
\begin{split}
a_{n,j,l}= & {1\over 2\pi i}
\sum_{k=1,k\ne j}^p\int_{C_k}{g_k(w)\over
(w-z_j)^{l+1}\prod_{s=1,s\ne k}^p(w-z_s)^{nm_s}}
{dw\over (w-z_k)^{nm_k+\rho_k}}+ \\ &
\int_{C_j}{g_j(w)\over\prod_{s=1,s\ne j}^p(w-z_s)^{nm_s}}
{dw\over (w-z_j)^{nm_j+\rho_j+l+1}},
\end{split}
\end{equation}
\begin{equation}
b_{n,j,l}= {1\over 2\pi i}
\sum_{k=q+1}^p\int_{C_k}
{\prod_{s=q+1,s\ne k}^p(w-z_s)^{nm_s}\over
(w-z_j)^{l+1}\prod_{s=1}^q(w-z_s)^{nm_s}}{g_k(w)dw\over (w-z_k)^{\rho_k-nm_k}},
\end{equation}
\begin{equation}
\begin{split}
c_{n,j,l}= &
\sum_{k=q+1,k\ne j}^p{1\over 2\pi i}\int_{C_k}
{\prod_{s=q+1,s\ne k}^p(w-z_s)^{(n+1)m_s}\over
(w-z_j)^{l+1}\prod_{s=1}^q(w-z_s)^{(n+1)m_s}}
{g_k(w)dw\over (w-z_k)^{\rho_k-(n+1)m_k}}+  \\ & {1\over 2\pi i}
\int_{C_j}{\prod_{s=q+1,s\ne j}^p(w-z_s)^{(n+1)m_s}\over
\prod_{s=1}^q(w-z_s)^{(n+1)m_s}}{g_j(w)
dw\over (w-z_j)^{\rho_j-(n+1)m_j+l+1}}.
\end{split}
\end{equation}
  From here, equations (\ref{coefvx}), (\ref{coefiiibis}) and
(\ref{coefiiibisbis})
  follow.
\end{proof}

\begin{remark}
\label{rem3}
     Let $z$ be a real or complex variable and
suppose that $(z-z_k)^{\rho_k}f(z)$
is $\rho_k-1-$times differentiable at $z_k$ for certain $\rho_k\in\Ns$.
Define
\begin{equation}
g(z)\equiv f(z)-\sum_{n=0}^{M}t_{n,m}^{(1)}(z)
{\prod_{k=1}^q(z-z_k)^{nm_k}\over
\prod_{k=q+1}^p(z-z_k)^{nm_k}}-
\sum_{n=0}^{M}t_{n,m}^{(2)}(z){\prod_{k=1}^q(z-z_k)^{(n+1)m_k}\over
\prod_{k=q+1}^p(z-z_k)^{(n+1)m_k}},
\end{equation}
where $M\equiv\lfloor$Max$\lbrace
(\rho_{q+1}-1)/m_{q+1},(\rho_{q+2}-1)/m_{q+2},\ldots,
(\rho_p-1)/m_p\rbrace\rfloor$ and
$t_{n,m}^{(1)}(z)$ and $t_{n,m}^{(2)}(z)$
are the polynomials defined in (\ref{polyii}), (\ref{polyiii}),
(\ref{coefiiibis}) and (\ref{coefiiibisbis}).
Then, the thesis of Proposition \ref{prop2}
holds for $f(z)$ replaced by $g(z)$.

Moreover, if
$\prod_{k=q+1}^p(z-z_k)^{\rho_k}f(z)$ is an
analytic function in $\Omega$, then the thesis of Theorem \ref{theo2}
applies to $g(z)$.
\end{remark}

\section{Discussion and concluding remarks}

In an earlier paper \cite{nicoii} we have discussed the theory
of two-point Taylor expansions, two-point Laurent
expansions and two-point Taylor-Laurent expansions. In
the present paper we have generalized these two-point
cases to multi-point cases. We have given details on the
regions of convergence and on representations of the
coefficients and the remainders of the expansions in terms
of Cauchy-type integrals.

Multi-point Taylor expansions are related with topics from
interpolation theory, in particular with the Newton
interpolation theory with applications in numerical
analysis. For example, applications can be found in
initial and boundary value problems in connection with
ordinary differential equations and in numerical
quadrature of integrals.

  From the point of view of interpolation theory detailed
information on multi-point expansions can be found in
\cite{walsh}, Chapters 3 and 8. The theory of several-point
Taylor expansions is discussed in Chapter 3 of \cite{walsh},
although in a setting that is different from our approach.
Our approach gives explicit Cauchy-type integrals of
coefficients and remainders which cannot be found in
Walsh's approach. In particular, we cannot find explicit
formulas for the polynomials $q_{n,m}(z)$ of formula (15)
as we have in (16)-(17). Knowledge of these explicit
formulas is necessary to construct asymptotic expansions
of integrals with several saddle points.

In addition to this, our Laurent and Taylor-Laurent
expansions are new. They have a formal similarity with
the rational approximations of Chapter 8 of \cite{walsh}: they
involve negative powers of $z$. But they are completely
different. The rational approximations, in particular the
Pad\'e-type approximations $P_n(z)/Q_m(z)$ are of
interpolatory type. These are generalizations of the
Taylor polynomial at several points: a quotient of
polynomials instead of a polynomial. However, our
expansions (21) or (32) have a different form and a
different approximation property: they approach not only
at regular points like Pad\'e-type approximations but also
at singular points of $f(z)$. And of course, the regions
and convergence properties in \cite{walsh} are different from
ours.

Apart from applying the present results in problems from
interpolation theory, in particular in problems from
numerical analysis, we expect to find applications in
asymptotic analysis of integrals, which application area
is our main motivation; see [5]. In that paper certain
orthogonal polynomials have been considered and we have
given new convergent expansions that also have an
asymptotic property for large values of a parameter (the
degree $n$ of the polynomials). Orthogonal polynomials
and special functions can be studied when the variable
and several parameters  are large. In that case more
than one or two so-called critical points occur that may
give the main contributions to the integral, and
expansions of analytic functions at these points gives
again the possibility of constructing new convergent
expansions with an asymptotic property. This method
avoids the complicated conformal mapping of the phase
function of the integral into a standard form (say a cubic
or higher polynomial). In addition, when the critical
points are multiple poles, Laurent-type expansions may be
considered. A few application areas are mentioned in the
Introduction, see the integral in (\ref{besel}), which we expect
to approximate in terms of Airy functions and the Pearcey
integral (\ref{pearcy}) and its derivative with respect to
$x$ and $y$.

\bibliographystyle{amsplain}

\end{document}